\documentclass[reqno,12pt]{article} \usepackage{fontenc}
\usepackage{amsfonts}  \usepackage{amssymb}

\newtheorem{thm}{Theorem}[section]

\newtheorem{rem}[thm]{Remark}
\newtheorem{ex}[thm]{Example}

\newtheorem{define}[thm]{Definition}
\usepackage{amssymb}
\usepackage{amsfonts}

\begin{document}


\begin{center}
{\bf \Large
Fuzzy Hom-Lie Subalgebras of Hom-Lie Algebras}\\~~\\
{Shadi Shaqaqha}\\~\\
Yarmouk University, Irbid, Jordan\\
shadi.s@yu.edu.jo
\end{center}

\vspace{.4cm} \begin{quote}
{\small \bf Abstract.} {\small In this paper, we introduce the concept of fuzzy Hom-Lie subalgebras (ideals) of Hom-Lie algebras and we investigate some of their properties. We study the relationship between fuzzy Hom-Lie subalgebras (resp. ideals) and Hom-Lie subalgebras (resp. ideals). For a finite number of fuzzy Hom-Lie subalgebras, we construct a new fuzzy hom-Lie subalgebras on their direct sum. Finally, The properties of fuzzy Hom-Lie subalgebras and fuzzy Hom-Lie ideals under morphisms of Hom-Lie algebras are studied.}
\end{quote}

\vspace{.4cm} \noindent {\bf Keywords:} Hom-Lie algebras; morphism of Hom-Lie algebras; direct sum, fuzzy set; fuzzy Hom-Lie subalgebra; fuzzy Hom-Lie ideal.

\vspace{.4cm} \noindent
\section{Introduction} \label{111}
The notion of Hom-Lie
	algebras was originally introduced by Hartwig, Larsson, and Silvestrov in 2006 \cite{HLA}. It is one of
	generalizations of the concept of classical Lie algebras. In recent years, they have become
	an interesting subject of mathematics and physics. We refer for more details on Hom-Lie algebras to \cite{subalgebra, Makhlouf2, Makhlouf3, Kdaisat, Shadi3}.\\
The idea of fuzzy sets was firstly introduced by Zadeh \cite{Zadeh}. A fuzzy set on a nonempty set $X$ is a map, called membership function, $\mu:X\rightarrow [0, 1]$. Note that in the the classical set theory we write $\mu(x)=1$ if $x\in X$, and $\mu(x)=0$ if $x\notin X$. Applications of the fuzzy set theory can be found in artificial intelligence, computer science,decision
theory, logic and management science, etc..\\
The study of fuzzy Lie subalgebras of Lie algebras was initiated by Yehia \cite{Yehia1} in 1996. Later fuzzy sets (and more generally intuitionistic fuzzy sets and complex fuuzy sets) have been applied in various directions in Lie algebras by many authors (see e.g. \cite{Akram,Akram1, Akram2, Davvaz, Shadi1, Shadi2}, and references therein) .\\
In this paper we describe fuzzy Hom-Lie algebras.
\section{Preliminaries} \label{222} 
Let $F$ be a ground field.	A Hom-Lie algebra over $F$ is a triple $(L, ~[~,~],~\alpha)$ where $L$ is a vector space over $F$, $\alpha: L\rightarrow L$ is a linear map, and $[~,~]: L\times L\rightarrow L$ is a bilinear map (called a bracket), satisfying the following properties:
	\begin{itemize}
		\item[(i)] $[x,~y]= -[y,~x]$ for all $x, y\in L$ (skew-symmetry property).
		\item[(ii)] $[\alpha(x), ~[y,~z]]+ [\alpha(y),~[z,~x]]+[\alpha(z),~[x,~y]]=0$ ,
		for all $x, y, z \in L$ (Hom-Jacobi identity).		
	\end{itemize}
It is clear that every Lie algebra is a Hom-Lie algebra by setting $\alpha= id_L$ (The identity map). For a Hom-Lie algebra $L$ over a field $F$ of characteristic $\neq2$, as in the setting of Lie algebras one can show that $[x,~x]=0$ for each $x\in L$. Also for an arbitrary Hom-Lie algebra $L$, we have $[x,~0]=[0,x]=0$ for each $x\in L$.
\begin{ex}
	Let $L$ be a vector space over $F$ and $[~,~]: L\times L\rightarrow L$ be any skew-symmetric bilinear map. If $\alpha: L\rightarrow L$ is the zero map, then  $(L, ~[~,~],~\alpha)$ is a Hom-Lie algebra.\hfill$\blacksquare$  	
\end{ex}
Let  $(L, ~[~,~],~\alpha)$ be a Hom-Lie algebra. A subspace $H$ of $L$ is a Hom-Lie subalgebra if $\alpha(H)\subseteq H$ and $[x,~y]\in H$ for all  $x,y\in H$.	A Hom-Lie subalgebra $H$ is  said to be a Hom-Lie ideal if $[x,~y]\in H$ for all  $x\in H$ and $y\in L$.\\
	Let $(L_1, ~[~,~]_1,~\alpha_1)$ and  $(L_2, ~[~,~]_2,~\alpha_2)$ be Hom-Lie algebras. A linear map $\varphi: L_1\rightarrow L_2$ is called a morphism of Hom-Lie algebras if the following two identities are satisfied:
	\begin{itemize}
		\item[(i)] $\varphi([x,~y]_1)= [\varphi(x),~\varphi(y)]_2$ for all $x, ~y\in L_1$.
		\item[(ii)]$\varphi\circ\alpha_1=\alpha_2\circ\varphi$.
	\end{itemize}  
Throughout this paper, $L$ is a Hom-Lie algebra over $F$.
\section{Fuzzy Hom-Lie Subalgebras and Fuzzy Hom-Lie Ideals} \label{333}
Let $a, b\in[0, 1]$. For the sake of simplicity we use the symbols $a\wedge b$ and $a\vee b$ to denote $\mathrm{min}\left\{a, b\right\}$ and $\mathrm{max}\left\{a, b\right\}$, respectively.
\begin{define}\label{CFL2}
A fuzzy set $\mu$ on $L$ is a {\em fuzzy Hom-Lie subalgebra} if the following conditions are satisfied for all $x, y\in L$, and $c\in F$:
\begin{itemize}
\item[(i)] $\mu(x+y)\geq \mu(x)\wedge \mu(y)$,
\item[(ii)] $\mu(c x)\geq \mu(x)$,
\item[(iii)] $\mu([x, y])\geq\mu(x)\wedge \mu(y)$,
\item[(iv)] $\mu(\alpha(x))\geq \mu(x)$.
\end{itemize}
\end{define}
If the condition (iii) is replaced by $\mu([x, y])\geq \mu(x) \vee\mu(y)$, then $\mu$ is called a {\em fuzzy Hom-Lie ideal} of $L$. Note that the condition $(ii)$ implies $\mu(x)\leq \mu(0)$ and $\mu(-x)\leq \mu(x)$ for all $x\in L$.\\
It is clear that if $\mu$ is a a fuzzy Hom-Lie ideal of $L$, then it is a fuzzy Hom-Lie subalgebra of $L$.\\
\begin{ex}\label{Homex2}
	Let $L$ be a vector space wiith basis $\{e_1,~e_2,~e_3\}$. We define the linear map $\alpha: L\rightarrow L$ by setting $\alpha(e_1)=e_2$ and $ \alpha(e_2) = \alpha(e_3)=0$. Let $[~,~]: L\times L\rightarrow L$ be the skew-symmetric bilinear map such that 
	$$[e_1,~e_2]=[e_2,~e_3]=0,~[e_1,~e_3]= e_1$$
	and $[e_i,~e_i]=0$ for all $i=1, 2, 3$. Then $(L, ~[~,~],~\alpha)$ is a Hom-Lie algebra. Indeed, for each $x, y\in L$, we have $[x,~y]$ is a scalar multiple of $e_1$. Also $\alpha(x)$ is a scalar multiple of $e_2$ for each $x\in L$. Therefore $[\alpha(x),~[y, ~z]]=0$ for each $x, y, z\in L$. This implies that the Hom-Jacobi identity is satisfied.\\
	We define $\mu$ as follows:

$$\mu(x)=\left\{
    \begin{array}{lr}
    0.8 &: x=0 \\
     0.4 &: x\in \mathrm{span}\{e_1, e_2\}-\{0\}\\
0.1 &: \mathrm{otherwise}.
 \end{array}
  \right.$$
Then $\mu$ is a fuzzy Hom-Lie ideal of $L$. \hfill$\blacksquare$	 
\end{ex}
\section{Relations Between Fuzzy Hom-Lie Ideals and Hom-Lie Ideals}
Let $V$ be a vector space and $\mu$ be a fuzzy set on it. For $t\in[0, 1]$ the set $U(\mu, t)= \left\{x\in V~|~\mu(x)\geq t\right\}$ is called an upper level of $\mu$. The following theorem will show a relation between fuzzy Hom-Lie subalgebras of $L$ and Hom-Lie subalgebras of $L$.
\begin{thm}\label{CFL3}
Let $\mu$ be a fuzzy subset of $L$. Then the following statements are equivalent:
\begin{itemize}
\item[(i)] $\mu$ is a fuzzy Hom-Lie subalgebra of $L$,
\item[(ii)] the non empty set $U(\mu, t)$ is a Hom-Lie subalgebra of $L$ for every $t\in \mathrm{Im}(\mu)$.
\end{itemize}
\end{thm}
{\it Proof.~}
Let $t\in\mathrm{Im}(\mu)$, and let $x, y\in U(\mu, t)$, and $c\in F$. As $\mu$ is a fuzzy Hom-Lie subalgebra of $L$, we have $\mu(x+y)\geq\mu(x)\wedge\mu(y)\geq t$, $\mu(c.x)\geq\mu(x)\geq t$, $\mu(\alpha(x))\geq \mu(x)\geq t$, and $\mu([x, y])\geq\mu(x)\wedge\mu(y)\geq t$, and so $x+y$, $\alpha x$, and $[x, y]$ are elements in $U(\mu, t)$. Conversely, let $U(\mu, t)$ be Hom-Lie subalgebras of $L$ for every $t\in \mathrm{Im}(\mu)$. Let $x, y\in L$ and $c\in F$. We may assume $\mu(y)\geq \mu(x)=t_1$, so $x, y\in U(\mu, t_1)$. As $U(\mu, t_1)$ is a subspace of $L$, we have $c. x$ and $x+y$ are in $U(\mu, t_1)$, and so $\mu(c. x)\geq t_1= \mu(x)$ and $\mu(x+y)\geq t_1 =\mu(x)\wedge\mu(y)$. Since $U(\mu, t_1)$ is a Hom-Lie subalgebra of $L$, we have $[x, y]$ and $\alpha(x)$ are in $U(\mu, t_1)$. Hence, $\mu([x, y])\geq t_1 =\mu(x)\wedge\mu(y)$, and $\mu(\alpha(x)\geq t_1=\mu(x)$.\hfill$\Box$
\begin{thm}\label{CFL4}
Let $\mu$ be a fuzzy subset of $L$. Then the following statements are equivalent:
\begin{itemize}
\item[(i)] $\mu$ is a fuzzy Hom-Lie ideal of $L$,
\item[(ii)] the non empty set $U(\mu, t)$ is a Hom-Lie ideal of $L$ for every $t\in \mathrm{Im}(\mu)$.
\end{itemize}
\end{thm}
{\it Proof.~} Let $\mu$ be a fuzzy Hom-Lie ideal of $L$. Then it is a fuzzy Hom-Lie subalgebra of $L$. According to the theorem above, every $x, y\in U(\mu, t)$ and $c\in F$ we have $x+y, c.x, \alpha(x)$ are in $U(\mu, t)$. For $x\in L$ and $y\in U(\mu, t)$, we find $\mu([x, y])\geq \mu(x)\vee \mu(y)\geq \mu(x)\geq t$. That is $[x, y]\in U(\mu, t)$. Conversely, assume that every $U(\mu, t)\neq \Phi$ is a Hom-Lie ideal of $L$, then $U(\mu, t)$ is a Hom-Lie subalgebra. Thus , we can proceed as in the theorem above and the only difference appears in the proof of the following statement:
$$\mu([x, y]\geq \mu(x)\vee \mu(y)~\forall x, y\in L.$$
Let $x, y\in L$. Without loss of generality, we may assume that $\mu(x)\geq \mu(y)$. Set $t_0=\mu(x)$. Hence $x\in U(\mu, t_0)$. As $U(\mu, t_0)$ is a Hom-Lie ideal of $L$, we have $[x, y]\in L$. This implies that $\mu([x, y])\geq \mu(x)\vee \mu(y)$.\hfill$\Box$

Let $V$ be a vector space. For $t\in [0, 1]$ and a fuzzy set $\mu$ on $V$, the set $U(\mu^>, t)= \left\{x\in V~|~\mu(x)> t\right\}$ is called a {\em strong upper level} of $\mu$. We have the following result.
\begin{thm}\label{CFL3B}
Let $\mu$ be a fuzzy subset of $L$. Then the following statements are equivalent:
\begin{itemize}
\item[(i)] $\mu$ is a fuzzy Hom-Lie subalgebra of $L$,
\item[(ii)] the strong upper level $U(\mu^>, t)$ is a subalgebra of $L$ for every $t\in \mathrm{Im}(\mu)$.
\end{itemize}
\end{thm}
{\it Proof.~}
For $t\in\mathrm{Im}(\mu)$, let $x, y\in U(\mu^>, t)$, and $c\in F$. As $\mu$ is a fuzzy Hom-Lie subalgebra of $L$, we have $\mu(x+y)\geq\mu(x)\wedge\mu(y)>t$, $\mu(c.x)>\mu(x)\geq t$, $\mu(\alpha(x))>\mu(x)$, and $\mu([x, y])\geq\mu(x)\wedge\mu(y)> t$. Consequently, $x+y$, $c.x$, $\alpha(x)$, and $[x, y]$ are elements in $U(\mu^>, t)$. Conversely, assume that for every $t\in \mathrm{Im}(\mu)$ we have $U(\mu^>, t)$ is a Hom-Lie subalgebra of $L$. Let $x, y\in L$ and $c\in F$. We need to show that the conditions of Definition \ref{CFL2} are satisfied. If $\mu(x)=0$ or $\mu(y)=0$, then $\mu(x+y)\geq 0=\mu(x)\wedge\mu(y)$. Suppose that $\mu(x)\neq 0$ and $\mu(y)\neq 0$. Suppose to the contrary that $\mu(x+y)$ and $\mu([x, y]$ are less than $\mu(x)\wedge \mu(y)$. Let $t_0$ be the greatest lower bound of the set $\{t~|~t<\mu(x)\wedge \mu(y)\}$. Since $x, y\in U(\mu^>, t_0)$, we have $x+y, [x, y]\in U(\mu^>, t_0)$, and hence $\mu(x+y), \mu([x, y])> t_0$. This contradicts that there is no element $a\in L$ with $t_0< \mu(a)< \mu(x)\wedge \mu(y)$. This shows that  $\mu(x+y), \mu([x, y])\geq\mu(x)\wedge \mu(y)$. Again let $t_0$ be the largest number  of $[0, 1]$ such that $t_0<\mu_A(x)$ and there is no $a\in L$ with $t_0<\mu_A(a)<\mu_A(x)$. As $U(\mu^>, t_0)$ is a Hom-Lie subalgerba, we have $c.x, \alpha(x)$ are in $U(\mu^>, t_0)$, and so $\mu(c.x)> t_0$ and $\mu(\alpha(x))> t_0$. Thus $\mu(c. x)$ and $\mu(\alpha(x))$  are greater than or equal to $\mu(x)$.\hfill$\Box$

Using almost the same argument one can show the following result.
\begin{thm}\label{CFL3C}
Let $\mu$ be a fuzzy subset of $L$. Then the following statements are equivalent:
\begin{itemize}
\item[(i)] $\mu$ is a fuzzy Hom-Lie ideal of $L$,
\item[(ii)] every strong upper level $U(\mu^>, t)$ is a Hom-Lie ideal of $L$ for every $t\in \mathrm{Im}(\mu)$.
\end{itemize}
\end{thm}
\section{Direct Sum of Fuzzy Hom-Lie Subalgebras}
	Given $n$ Hom-Lie algebras $(L_i, ~[~,~]_i,~\alpha_i)$ ,$i=1, \ldots,n$, then $$(L_1\oplus L_2\oplus \ldots\oplus L_n,~[~, ~],~\alpha_1+ \alpha_2+ \ldots+\alpha_n)$$
	is a Hom-Lie algebra by setting
	$$[~, ~]~:~(L_1\oplus L_2\oplus \ldots\oplus L_n)\times( L_1\oplus L_2\oplus \ldots\oplus L_n)\rightarrow (L_1\oplus L_2\oplus \ldots\oplus L_n)$$
	$$((x_1,~x_2, \ldots,~x_n),~(y_1,~y_2, \ldots,~y_n))\mapsto ([x_1,~y_2]_1,~[x_1,~y_2]_2, \ldots,~[x_n,~y_n]_n),$$ 
	and the linear map $$(\alpha_1+ \alpha_2+ \ldots+\alpha_n)~:~(L_1\oplus L_2\oplus \ldots\oplus L_n)\rightarrow(L_1\oplus L_2\oplus \ldots\oplus L_n)$$
	$$(x_1,~x_2, \ldots,~x_n) \mapsto (\alpha_1(x_1),~\alpha_2(x_2), \ldots,~\alpha_n(x_n)).$$ 
In the special case where $n=2$, we obtain \cite[Proposition 2.2]{directsum} (see \cite{Kdaisat}).\\
Let $(L_1, [~,~]_1, \alpha_1), (L_2, [~,~]_2, \alpha_2), \ldots, (L_n, [~,~]_n, \alpha_n)$ be Hom-Lie algebras. Suppose that $\mu_1, \mu_2, \ldots, \mu_n$ are fuzzy subsets of $L_1, L_2, \ldots, L_n$, respectively. Then the generalized Cartesian sum of fuzzy sets induced by $\mu_1, \mu_2, \ldots, \mu_n$ on $L_1\oplus L_2\oplus \cdots \oplus L_n$ is 
$$\mu_1\oplus \mu_2\oplus \cdots\oplus \mu_n:L_1\oplus L_2\cdots \oplus L_2\rightarrow [0,1];~(x_1, x_2, \ldots, x_n)\mapsto \mu_1(x_1)\wedge \mu_2(x_2)\wedge \mu_n(x_n).$$
\begin{thm}\label{proifnls1}
	Let $(L_1, [~,~]_1, \alpha_1), (L_2, [~,~]_2, \alpha_2), \ldots, (L_n, [~, ~]_n, \alpha_n)$ be Hom-Lie algebras. Let $\mu_1, \mu_2, \ldots, \mu_n$ be fuzzy Hom-Lie aubalgebras of $L_1, L_2, \ldots, L_n$, respectively. Then $\mu_1\oplus \mu_2\oplus\cdots \oplus \mu_n$ is a fuzzy Hom-Lie subalgebra of $L_1\oplus L_2\oplus \cdots \oplus L_n$.
\end{thm}
{\it Proof.~} Let $(x_1, x_2, \ldots, x_n), (y_1, y_2, \ldots, y_n)\in L_1\oplus L_2\oplus \cdots \oplus L_n$. Then
\begin{eqnarray}
	\mu_1([(x_1, x_2, \ldots, x_n), (y_1, y_2, \ldots, y_n)])&=& (\mu_1\oplus \mu_2\oplus\cdots \oplus \mu_n)([x_1, y_1]_1, [x_2, y_2]_2, \ldots, [x_n, y_n]_n)\nonumber\\
	&=&\mu_1([x_1, y_1]_1)\wedge \mu_2([x_2, y_2]_2)\wedge \cdots\wedge \mu_n([x_n, y_n]_n)\nonumber\\
	&\geq& (\mu_1(x_1)\wedge\mu_1(y_1)\wedge \mu_2(x_2)\wedge \mu_2(y_2)\cdots \wedge \mu_n(x_n)\wedge \mu_n(y_n)\nonumber \\
	&=&(\mu_1\oplus \cdots\oplus \mu_n)((x_1, \ldots, x_2))\wedge (\mu_1\oplus \cdots\oplus \mu_n)((y_1, \ldots, y_n)).\nonumber
\end{eqnarray} 
Also, 
\begin{eqnarray}
  (\mu_1\oplus \cdots \oplus \mu_n)(\alpha_1+\alpha_2+\cdots +\alpha_n)(x_1, x_2, \ldots, x_n)&=&(\mu_1\oplus \cdots \oplus \mu_n)(\alpha_1(x_1), \alpha_2(x_2), \ldots, \alpha_n(x_n))\nonumber\\
  &=&\nonumber\mu_1(\alpha_1(x_1))\wedge \mu_2(\alpha_2(x_2))\wedge \cdots \wedge \mu_n(\alpha_n(x_n))\nonumber\\
  &\geq& \mu_1(x_1)\wedge \mu_2(x_2)\wedge \cdots \wedge \mu_n(x_n)\nonumber\\
  &=&(\mu_1\oplus \mu_2\oplus \cdots \mu_n)*x_1, x_2, \ldots, x_n).\nonumber
\end{eqnarray}
The rest of the proof is similar to the proof of \cite[Theorem 5.2]{Huang}, so we omit it.
\hfill $\Box$

However the direct sum of fuzzy-Hom Lie ideals of Hom-Lie algebras $L_1$ and $L_2$ is not nesaccary to be a fuzzy Hom-Lie ideal of the Hom-Lie algebra $L_1\oplus L_2$.
\begin{rem}
In \cite{Shadi7}, we introduced and studied infinite direct product of Hom-Lie algebras. One can consider fuzzy Hom-Lie subalgebras of such hom-Lie algebras.  
\end{rem}
\section{On Fuzzy Hom-Lie algebras and Hom-Lie Algebras Morphisms}
Suppose $f: X\rightarrow Y$ is a function. If $\mu_B$ is a fuzzy set of $Y$, then we can define a fuzzy set on $X$ induced by $f$ and $\mu_A$ by setting $\mu_{f^{-1}(B)}(x)=\mu_B(f(x))$ for any $x\in X$. Also if $\mu_A$ is a fuzzy set on $X$, then 
$$\mu_{f(A)}(y)=\left\{
    \begin{array}{lr}
    \mathrm{sup}_{x\in f^{-1}(y)}\left\{\mu_A(x)\right\} &: y\in f(X) \\
     0 &: y\notin f(X)
 \end{array}
  \right.$$
is a fuzzy set on $Y$ induced by $f$ and $\mu_A$ (See for example \cite{Shadi3}). The following theorem was obtained by Kim and Lee in \cite{Kim} in the setting of Lie algebras. We extend it to Hom-Lie algebra case.
\begin{thm}\label{CFL9}
Let $f: (L_1, [~,~]_1, \alpha_1)\rightarrow (L_2, [~,~]_2, \alpha_1)$ be a morphism of Hom-Lie algebras. If $B=\mu_B$ is a fuzzy Hom-Lie subalgebra (resp. ideal) of $L_2$, then the fuzzy set $f^{-1}(B)$ is also a fuzzy Hom-Lie subalgebra (resp. ideal) of $L_1$.
\end{thm}
{\it Proof.~}
Let $x_1, x_2\in L_1$. Then
\begin{eqnarray*}
\mu_{f^{-1}(B)}(x_1+x_2)&=&\mu_B(f(x_1+x_2))\\
&=& \mu_B(f(x_1)+f(x_2))~~~(f~\mathrm{is~linear})\\
&\geq & \mu_B(f(x_1))\wedge \mu_B(f(x_2))~(\mu_B~\mathrm{is ~a~ fuzzy ~Hom-Lie~subalgebra})\\
&=&\mu_{f^{-1}(B)}(x_1)\wedge \mu_{f^{-1}(B)}(x_2),
\end{eqnarray*}
and
\begin{eqnarray*}
\mu_{f^{-1}(B)}([x_1,x_2])&=&\mu_B(f([x_1, x_2])\\
&=& \mu_B([f(x_1), f(x_2)])~~~(f~\mathrm{is~morphism})\\
&\geq & \mu_B(f(x_1))\wedge \mu_B(f(x_2))~(\mu_B~\mathrm{is~a~fuzzy~Hom-Lie~subalgebra})\\
&=&\mu_{f^{-1}(B)}(x_1)\wedge \mu_{f^{-1}(B)}(x_2).
\end{eqnarray*}
Let $x\in L_1$ and $c\in F$. Then 
\begin{eqnarray*}
\mu_{f^{-1}(B)}(c .x)&=&\mu_B(f(c.x))\\
&=& \mu_B(c.f(x))~~~(f~\mathrm{is~linear})\\
&\geq & \mu_B(f(x))~(\mu_B~\mathrm{is ~a~fuzzy~Hom-Lie~subalgebra})\\
&=&\mu_{f^{-1}(B)}(x),
\end{eqnarray*}
and
\begin{eqnarray*}
\mu_{f^{-1}(B)}(\alpha_1(x))&=&\mu_B(f(\alpha_1(x)))\\
&=& \mu_B(\alpha_2(f(x)))~~~(f~\mathrm{is~a~morphism~of~Hom-Lie~algebras})\\
&\geq & \mu_B(f(x))~(\mu_B~\mathrm{is ~a~fuzzy~Hom-Lie~subalgebra})\\
&=&\mu_{f^{-1}(B)}(x),
\end{eqnarray*}
The case of fuzzy Hom-Lie ideal is similar to show.\hfill$\Box$ 

If $f:L_1\rightarrow L_2$ is a Lie algebra homomorphism and $A=\mu_A$ is a fuzzy subalgebra of $L_1$, then the image of $A$, $f(A)$ is a fuzzy subalgebra of $f(L_1)$ (\cite{Kim}). In the following theorem we establish an analogue result for the case of Hom-Lie algebras.

\begin{thm}\label{CFL10}
Let $f:(L_1, [~,~]_1, \alpha_1)\rightarrow (L_2, [~,~]_2, \alpha_2)$ be a morphism from $L_1$ onto $L_2$. If $A=\mu_A$ is a fuzzy Hom-Lie subalgebra of $L_1$, then $f(A)$ is also a fuzzy Hom-Lie subalgebra of $L_2$.
\end{thm}
{\it Proof.~} Let $y_1, y_2\in L_2$. As $f$ is onto, there are $x_1, x_2\in L_1$ such that $f(x_1)=y_1$ and $f(x_2)=y_2$. We have 
$$\{x_1+x_2~|~x_1\in f^{-1}(y_1)~\mathrm{and}~x_2\in f^{-1}(y_2)\}\subseteq \{x~|~x\in f^{-1}(y_1+y_2)\},$$
and 
$$\{[x_1, x_2]_1~|~x_1\in f^{-1}(y_1)~\mathrm{and}~x_2\in f^{-1}(y_2)\}\subseteq \{x~|~x\in f^{-1}([y_1, y_2]_2)\}.$$
Now, we find
\begin{eqnarray}
  \mu_{f(A)}(y_1+y_2)&=&\mathrm{sup}_{x\in f^{-1}(y_1+y_2)}\{\mu_A(x)\}\nonumber\\
  &\geq& \{\mu_A(x_1+x_2)~|~x_1\in f^{-1}(y_1)~\mathrm{and}~x_2\in f^{-1}(y_2)\}\nonumber\\
  &\geq &\mathrm{sup}\{\mu_A(x_1)\wedge \mu_A(x_2)~|~x_1\in f^{-1}(y_1)~\mathrm{and}~x_2\in f^{-1}(y_2)\}\nonumber\\
  &=& \mathrm{sup}_{x_1\in f^{-1}(y_1)}\{\mu_A(x_1)\}\wedge \mathrm{sup}_{x_2\in f^{-1}(y_2)}\{\mu_A(x_2)\}\nonumber\\
  &=&\mu_{f(A)}(y_1)\wedge \mu_{f(A)}(y_2).\nonumber
\end{eqnarray}
Also, 
\begin{eqnarray}
  \mu_{f(A)}([y_1, y_2]_2)&=&\mathrm{sup}_{x\in f^{-1}([y_1, y_2]_2)}\{\mu_A(x)\}\nonumber\\
  &\geq &\{\mu_A([x_1, x_2]_1)~|~x_1\in f^{-1}(y_1)~\mathrm{and}~x_2\in f^{-1}(y_2)\}\nonumber\\
  &\geq&\mathrm{sup}\{\mu_A(x_1)\wedge \mu_A(x_2)~|~x_1\in f^{-1}(y_1)~\mathrm{and}~x_2\in f^{-1}(y_2)\}\nonumber\\
  &=& \mathrm{sup}_{x_1\in f^{-1}(y_1)}\{\mu_A(x_1)\}\wedge \mathrm{sup}_{x_2\in f^{-1}(y_2)}\{\mu_A(x_2)\}\nonumber\\
  &=&\mu_{f(A)}(y_1)\wedge \mu_{f(A)}(y_2).\nonumber
\end{eqnarray}
For $y\in L_2$ and $c\in F$, we find
$$\{c.x~|~x\in f^{-1}(y)\}\subseteq \{x~|~x\in f^{-1}(c.y)\},$$
and
$$\{\alpha_1(x)~|~x\in f^{-1}(y)\}\subseteq \{x~|~x\in f^{-1}(\alpha_2(y))\}.$$
and so

\begin{eqnarray}
  \mu_{f(A)}(c.y)&=&\mathrm{sup}_{x\in f^{-1}(c.y)}\{\mu_A(x)\}\nonumber\\
  &\geq& \{\mu_A(c.x)~|~x\in f^{-1}(c.y)\}\nonumber\\
  &\geq& \{\mu_A(x)~|~x\in f^{-1}(y)\}\nonumber\\
  &=&\mu_{f(A)}(y),\nonumber
\end{eqnarray}
also,
\begin{eqnarray}
  \mu_{f(A)}(\alpha_2(y))&=&\mathrm{sup}_{x\in f^{-1}(\alpha_2(y))}\{\mu_A(x)\}\nonumber\\
  &\geq&\{\mu_A(\alpha_1(x))~|~x\in f^{-1}(\alpha_2(y))\}\nonumber\\
  &\geq& \{\mu_A(x)~| x\in f^{-1}(y)\}\nonumber\\
  &=&\mu_{f(A)}(y).\nonumber
\end{eqnarray}
\hfill$\Box$

Chung-Gook Kim and Dong-Soo Lee (\cite{Kim}) proved if $\varphi:L\rightarrow L'$ is a surjective Lie algebra homomorphism and $A=\mu_A$ is a fuzzy ideal of $L$, then $\varphi(A)$ is a fuzzy ideal of $L'$. We will extend the result to fuzzy Hom-Lie algebra case.

\begin{thm}\label{CFL11}
Let $f:(L_1, [~,~]_1, \alpha_1)\rightarrow (L_2, [~,~]_2, \alpha_2)$ be an onto morphism of Hom-Lie algebras. If $A=\mu_A$ is a fuzzy Hom-Lie ideal of $L_1$, then $f(A)$ is also a fuzzy Hom-Lie ideal of $L_2$.
\end{thm}
{\it Proof.~}
The proof is similar to the proof of the theorem above. We only need to show that $\mu_{f(A)}([y_1, y_2]_2)\geq \mu_{f(A)}(y_1)\vee \mu_{f(A)}(y_2)$ for all $y_1, y_2\in L_2$. Let $y_1, y_2\in L_2$, and assume, by contradiction, that $\mu_{f(A)}([y_1, y_2]_2)<\mu_{f(A)}(y_1)\vee \mu_{f(A)}(y_2)$. Then $\mu_{f(A)}([y_1, y_2]_2)< \mu_{f(A)}(y_1)$ or $\mu_{f(A)}([y_1, y_2]_2)< \mu_{f(A)}(y_2)$. We may assume, without loss of generality, that $\mu_{f(A)}([y_1, y_2]_2)< \mu_{f(A)}(y_1)$. Choose a number $t\in [0, 1]$ such that $\mu_{f(A)}([y_1, y_2]_2)<t<\mu_{f(A)}(y_2)$. There is $a\in f^{-1}(y_1)$ with $\mu_A(a)>t$. As $f$ is onto, there exists $b\in f^{-1}(y_2)$. We note that
$$f([a, b]_1)=[f(a), f(b)]_2=[y_1, y_2]_2.$$
Thus,
$$\mu_{f(A)}([y_1, y_2]_2)\geq \mu_A([a, b]_1)\geq \mu_A(a)\vee \mu_A(b)>t>\mu_{f(A)}([y_1, y_1]_2).$$
Contradiction.\hfill$\Box$


\bigskip \bigskip

\begin{thebibliography}{99}

\bibitem{Akram book} {\bf M. Akram}, {\em Fuzzy Lie algebras},	Springer Nature Singapore Pte Ltd. (2018).	
\bibitem{Akram} {\bf M. Akram}, {\em Intuitionistic fuzzy Lie algebras},  Southeast Asian Bulletin of Mathematics {\bf 31} (2007), 843-855.
\bibitem{Akram1} {\bf M. Akram}, {\em Intuitionistic (S,T)-fuzzy Lie ideals of Lie algebras}, Quasigroups Relat. Systems {\bf 15}
(2007), 201-218.
\bibitem{Akram2} {\bf M. Akram}, {\em Intuitionistic fuzzy Lie ideals of Lie algebras}, Int. Journal of Fuzzy Math., {\bf 6} (2008), no. $4$, 991-1008.
\bibitem{subalgebra} {\bf J. Casas, M. Insua and N. Pacheco}, {\em On universal central extensions of Hom-Lie algebras}, Hacettepe Journal of Mathematics and Statistics  {\bf 44} (2015), no. $2$, 277-288.
\bibitem{Davvaz} {\bf B. Davvaz and WA. Dudek}, {\em Fuzzy $n$-Lie algebras},  J Generalized Lie Theory Appl {\bf 11} (2017), 1-6.
\bibitem{HLA} {\bf J. Hartwig, D. Larsson,  and S. Silvestrov}, {\em Deformations of Lie algebras using $\sigma$-derivations}, Journal of Algebra {\bf 295} (2006), no. $2$, 314-361.
\bibitem{Huang} {\bf C. E. Huang and F. G. Shi}, \textit{On the fuzzy dimensions of fuzzy vector spaces}, Iranian Journal of Fuzzy Systems, {\bf 9} (4), 141–150.
\bibitem{Kdaisat} {\bf N. Kdaisat}, {\em On Hom-Lie algebras}, Master thesis, Yarmouk University (2021).
\bibitem{Kim} {\bf Chung-Gook Kim, Dong-Soo Lee}, {\em Fuzzy Lie ideals and fuzzy Lie subalgebras}, Fuzzy Sets and Systems {\bf 94} (1998), 101-104.
\bibitem{Makhlouf2} {\bf A. Makhlouf, S. Silvestrov}, {\em Notes on formal deformations of hom-associative and hom-Lie algebras}, Forum Math. {\bf } (2010), no. $4$, 715--739.
\bibitem{Makhlouf3} {\bf A. Makhlouf, S. Silvestrov}, {\em Hom-algebra structures}, J. Gen. Lie Theory Appl.{\bf 2} (2008), no. $2$, 51--64.

\bibitem{Shadi1} {\bf S. Shaqaqha}, {\em Complex fuzzy Lie algebras}, Jordan Journal of Mathematics and Statistics (JJMS) {\bf 13} (2020), no. $2$, 231 - 247.
\bibitem{Shadi2} {\bf S. Shaqaqha}, {\em On fuzzification of-$n$-Lie algebras}, Jordan Journal of Mathematics and Statistics, to appear.
\bibitem{Shadi3} {\bf S. Shaqaqha}, {\em Restricted Hom-Lie superalgebras}, Jordan Journal of Mathematics and Statistics (JJMS) {\bf 12} (2019), no. $2$, 233-255.
\bibitem{Shadi7} {\bf S. Shaqaqha and N. Kdaisat}, {\em More properties of (multiplicative)-Hom-Lie algebras}, submitted, 2022.
\bibitem{directsum}Y. Sheng, Representations of Hom-Lie algebras, {\it Algebras and Representation Theory}, {\bf 15}(6) (2012), 1081-1098.
\bibitem{Yehia1} {\bf Samy El-Badawy Yehia}, {\em Fuzzy ideals and fuzzy subalgebras of Lie algebras}, Fuzzy Sets and Systems {\bf 80} (1996), 237-244.
\bibitem{Yehia2} {\bf Samy El-Badawy Yehia}, {\em The adjoint representation of fuzzy Lie algebras}, Fuzzy Sets and Systems {\bf 119} (2001), 409-417.
\bibitem{Zadeh} {\bf L.A. Zadeh}, {\em Fuzzy sets}, Inform. Control {\bf 8} (1965), 338-358.

\end{thebibliography}
\end{document}